\documentstyle[twocolumn,aps,floats]{revtex}

\begin{document}

\draft \tolerance = 10000

\setcounter{topnumber}{1}
\renewcommand{\topfraction}{0.9}
\renewcommand{\textfraction}{0.1}
\renewcommand{\floatpagefraction}{0.9}

\twocolumn[\hsize\textwidth\columnwidth\hsize\csname
@twocolumnfalse\endcsname

\title{Generalized Riemann - Liouville fractional derivatives for multifractal sets}
\author{L.Ya.Kobelev   \\
Department of  Physics, Urals State University \\ Lenina Ave., 51,
Ekaterinburg 620083, Russia \\E-mail: leonid.kobelev@usu.ru  }
\maketitle

\begin{abstract}
The Riemann-Liouville fractional integrals and derivatives are
generalized for cases when fractional exponent $d$ are functions
of space and times coordinates (i.e. $d=d({r}(t),t)$).
\end{abstract}

 \vspace{1cm}

]

\section{Introduction}

Fractional derivatives and integrals (left-sided and right-sided)
Riemann - Liouville (see \cite{1}-\cite{3}) from functions $f(t)$
(defined on a class of generalized functions) are
\begin{equation} \label{eq1}
D_{+,t}^{d}f(t)=\frac{1}{\Gamma(n-d)}\left(
\frac{d}{dt}\right)^{n}\int_{a}^{t}
\frac{f(t^{\prime})dt^{\prime}}{(t-t^{\prime})^{d-n+1}}
\end{equation}
\begin{equation} \label{eq2}
D_{-,t}^{d}f(t)= \frac{(-1)^{n}}{\Gamma(n-d)}\left(
\frac{d}{dt}\right)
^{n}\int_{t}^{b}\frac{f(t^{\prime})dt^{\prime}}{(t^{\prime}-t)^{d-n+1}}
\end{equation}
where $\Gamma(x)$ is Euler's gamma function, and $a$ and $b$ are
some constants from $[0,\infty) $. In these definitions, as
usually, $n=\{d\}+1$ , where $\{d\}$ is the integer part of $d$ if
$d\geq 0$ (i.e. $n-1\le d<n$) and $n=0$ for $d<0$. Fractional
derivatives and the integrals (\ref{eq1})-(\ref{eq2}) allow to
use, instead of usual derivatives and integrals, the integral
functionals defined on a wide class of generalized functions. It
is very useful for the solution of a series of problems describing
stochastic and chaos processes, abnormal diffusion, quantum
theories of a field etc. \cite{4}-\cite{10}. It is possible to
consider appearance of integral in (\ref{eq1})-(\ref{eq2}), from
the physical point of view, as the result of taking into account
influence of the contributions from some physical processes
(characterized by the kernel $(t-t')^{-d+n-1}\Gamma^{-1}(n-d)$) in
earlier (left-side derivative) or later (right-hand derivative)
times, on function $f(t)$ that is, as the partial taking into
account the system memory about past or future times. The value of
fractional exponent $d$ characterizes the degree of the memory.
Let's consider multifractal set (without self-similarity) $S_{t}$
consisting from infinite number of subsets $s_{i}(t_{i})$, also
being multifractal. Each subset $s_{i}(t_{i})$ is compared with
fractional value (or number of values), describing its fractal
(fractional) dimension (box dimension, Hausdorff \cite{11} or
Renie \cite{12} dimension etc. - see, for example, \cite{13}),
depending from the numbers of a subset $s_{i}(t)$. Let the carrier
of measure of multifractal set $S_{t}$ be the set $R^{n}$. For
exposition of changes of a continuous function $f(t)$ defined on
subsets $s_{i}(t_{i})$ of set $S_{t}$, it is impossible to use
ordinary derivatives or Riemann - Liouville fractional derivatives
(\ref{eq1}), as the fractional dimension of sets $d$ on which
$f(t)$ is defined depends on $t_{i}$, that is on the choice of the
subset $s_{i}(t_{i)}$. There is a problem: how can the definition
(\ref{eq1})-(\ref{eq2}) be changed to feature small (or major)
changes of function $f(t)$ defined on sets $s_{i}(t_{i)}$? The
purpose of this paper is to present the generalization of the
Riemann - Liouville fractional derivatives (\ref{eq1})-(\ref{eq2})
in order to ajust them for functions defined on multifractal sets
with fractal dimension (fractional dimension) depending on the
coordinates and time.

\section{Generalized fractional derivatives and integrals}

We shall treat subsets $s_{i}(t_{i)}$ as the "points" $t_{i}$
(with a continuous distribution for different multifractal subsets
$s_{i}(t_{i)}$) of multifractal set $S_{t}$). Assume that the
function $d(t_{i})=d(t)$, describing their fractional dimension
(in some cases coinciding with local fractal dimension) as
function $t$ is continuous. For the elementary generalization
(\ref{eq1})-(\ref{eq2}) is used physical reasons and variable $t$
is interpreted as a time. For continuous functions $f(t)$
(generalized functions defined on the class of finitary functions
(see \cite{3})), the Riemann - Liouville fractional derivatives
also are continuous. So for infinitesimal intervals of time and
the functionals (\ref{eq1})-(\ref{eq2}) will vary on infinitesimal
quantity. For continuous function $d(t)$ the changes thus also
will be infinitesimal. It allows, as the elementary generalization
(\ref{eq1}) suitable for describing of changes $f(t)$ defined on
multifractal subsets $s(t)$, as well as in
(\ref{eq1})-(\ref{eq2}), to summate influence of a kernel of
integral $(t-t')^{-d(t)-n+1}\Gamma^{-1}(n-d(t))$, depending on
$d(t)$, in all points of integration and, instead of
(\ref{eq1})-(\ref{eq2}), introduce the following definitions
(generalized fractional derivatives and integrals (GFD)), taking
into account also the $d(t)$ dependence on time and vector
parameter ${\bf r}(t)$ (i.e. $d_t\equiv d_t({\bf r}(t),t)$)
\begin{equation}\label{eq3}
D_{+,t}^{d_t}f(t)=\left({\frac{d}{dt}}\right)^n \int\limits_a^t
{dt'\frac{f(t')}{\Gamma(n-d_t(t'))(t-t')^{d_t(t')-n+1}}}
\end{equation}
\begin{eqnarray}\label{eq4} \nonumber
& & D_{-,t}^{d_t}f(t)=(-1)^n \times \\ &\times& \left({\frac{d}{{dt}}}
\right)^n\int\limits_t^b {dt'\frac{{f(t')}}
{{\Gamma(n-d_t(t'))(t'-t)^{d_t(t')-n+1}}}}
\end{eqnarray}
In (\ref{eq3})-(\ref{eq4}), as well as in (\ref{eq1})-(\ref{eq2}),
$a$ and $b$ stationary values defined on an infinite axis (from
$-\infty$ to $\infty$), $a<b$ , $n-1 \leq d_{t}<n$,
$n=\{d_{t}\}+1$, $\{d_{t}\}$- the whole part of $d_{t}\geq 0$,
$n=0$ for $d_{t}<0$. The sole difference (\ref{eq3})-(\ref{eq4})
from (\ref{eq1})-(\ref{eq2}) is: $d_{t}=d_{t}({\bf r}(t),t)$-
fractional dimension (further will be used for it terms " fractal
dimension " (FD) or " the local fractal dimension (LFD) ") is the
function of time and coordinates, instead of stationary values in
(\ref{eq1})-(\ref{eq2}).

Similar to (\ref{eq3})-(\ref{eq4}), it is possible to define the
GFD, (coinciding for integer values of fractional dimension
$d_{\vec{\bf r}}({\bf r},t)$ with derivatives with respect to
vector variable ${\bf r}$) $D^{d_{{bf r}}}_{+,{\bf r}}f({\bf
r},t)$ respect to vector ${\bf r}(t)$ variables (spatial
coordinates). We pay attention, that definitions
(\ref{eq3})-(\ref{eq4})are a special case of Hadamard derivatives
\cite{14}.

\section{Fractional derivatives for $d({\bf r}(t),t) \to 1$}

For FD which have very small differences  from of integer values
it is possible approximate to change the GFD by the usual
derivatives and integrals. For an establishment of connection of
GFD with orderly derivatives we shall see (\ref{eq3}), for
example, for a case $d({\bf r}(t),t)=1+\varepsilon({\bf
r}(t),t)$,$\varepsilon \ll 1$, $d<1$, (if utilize the theorem of
the mean value of integral) as
\begin{eqnarray}\label{eq5} \nonumber
&& D_{+,t}^{1-\varepsilon}f(t)=\frac{\partial}{\partial t}\int\limits_0^t
{\frac{f(t-\tau)d\tau}{\Gamma(\varepsilon(t-\tau))(\tau\pm i\xi)^
{1-\varepsilon(t-\tau)}}}= \\ &=& \frac{\partial}{\partial t}[\tilde
f(t-\tau_{cp}
(t))\int\limits_0^t{\frac{d\tau}{(\tau)^{1-\varepsilon(t-\tau)}}}]
\end{eqnarray}
where $\tilde f=\Gamma^{-1}f$ and $t_{med}$- some value of $\tau$. As
$\varepsilon\to 0$ it is possible to estimate values of integral in
(\ref{eq5}) for minimum and maximal values of $\varepsilon$
($\varepsilon>0$)
\begin{equation}\label{eq6}
\int\limits_0^t{\frac{d\tau}{\tau^{1-\varepsilon_{\min}(t-\tau)}}}=
\frac{t^{\varepsilon_{\min}}}{\varepsilon_{\min}},\quad\int\limits_0^t
{\frac{d\tau}{\tau^{1-\varepsilon_{\max}(t - \tau)}}}=\frac{t^{\varepsilon
_{\max}}}{\varepsilon_{\max}}
\end{equation}
For selection from integrals (\ref{eq6}) the trms which are
independent from $\varepsilon$ (because of $\tilde{f}\sim
\varepsilon$) we use decomposition $t=1+\ln t+...$ We obtain
\begin{equation}\label{eq7}
D_{+,t}^{1-\varepsilon}f(t)\approx\frac{\partial f(t)}{\partial t}+
\frac{\partial\tilde f(t-\tau_{cp}(t)}{\partial t}\ln t+\frac{\tilde f(t-
\tau_{cp}(t)}{t}
\end{equation}
For major times $t=t_0+(t-t_0)$,$t-t_0\ll t_0$ the approximate
representation GFD (\ref{eq7}) through usual derivatives will accept a
view (if neglect by additions of order $\tilde f/t_0$,$(t-t_0)/t_0$, to
use the designation $\ln t_0 =\alpha$ and if it is accounted, that
$\tau_{cp}\ll t$ because of the basic contribution to integral (\ref{eq5})
is stipulated by small $\tau$)
\begin{equation}\label{eq8}
D_{+,t}^{1-\varepsilon}f(t)\approx\frac{\partial f(t)}{\partial t}+
\frac{\partial\alpha\tilde f(t)}{\partial t}
\end{equation}
In (\ref{eq8}) $\alpha$ play a role of a stationary value of parameter of
regularization and if change $\varepsilon$ (in $\tilde{f}=\varepsilon f$)
on quantity $\varepsilon\to\varepsilon\alpha^{-1}$, GFD (\ref{eq8}) is not
depends practically on this parameter.

We shall give below, another method of deduction  the relation
(\ref{eq8}) using an expansion $\tau^{1-\varepsilon}$ in a power
series of $\epsilon$ under sign of integral and again for
$d(r(t),t)$ with poorly difference from the whole value (but
$d>1$). Let's fractional dimension $d({\bf r}(t),t)$ is equal
unity with small value $\varepsilon$ ($d({\bf
r}(t),t)=1+\varepsilon({\bf r}(t),t)$,$\varepsilon\ll 1$) and
expand FD in (\ref{eq3}) in a power series on $\varepsilon$ by a
rule $(t-\tau)^{-\varepsilon}=1 \varepsilon\ln(t-\tau)+...$.
Restricted expansion FD by the first two members of a series, we
obtain for left-side fractional derivative (for $a=0$)
\begin{eqnarray}\label{eq9}
&& \nonumber D_{ + ,t}^{1 + \varepsilon } f(t) = \frac{{\partial ^2 }}
{{\partial t^2 }}\int\limits_0^t {\frac{{f(\tau )}} {{\Gamma (1 -
\varepsilon )(t - \tau )^{\varepsilon (r(\tau ),\tau )} }}d\tau }  \approx
\nonumber \\  \nonumber & \approx & \frac{\partial } {{\partial
t}}(\frac{1} {{\Gamma (1 - \varepsilon )}}f(t)) - \frac{{\partial ^2 }}
{{\partial ^{} t^2 }}\int\limits_0^t {\frac{{\varepsilon _{} f(\tau )d\tau
}} {{\Gamma (1 - \varepsilon )[(t - \tau )^{}  \pm i\varsigma ]}}}, \\ &&
\quad\varsigma \to 0
\end{eqnarray}
Integral in (\ref{eq9}) is considered as a generalized function with
determined regularizasion and after regularization of integral the
parameter of regularization $\alpha$ is picked by a requirement of the
best coincidence of approximate and exact results of integral calculation
the members of a first order at $\varepsilon$ (it is necessary to take its
real part, the parameter $\varsigma$ is necessary to put zero after
calculations). After an integration by parts and also using of a relation
$1/x   = P(1/x)$ or another regularization's we shall receive (if take
into account that integrals (\ref{eq1})-(\ref{eq2}) are real values and
define the fractal addendum to derivatives as  a coefficients at the
imaginary  parts of integrals)
\begin{equation}\label{eq10}
D_{+,t}^{1+\varepsilon}f(t)=\frac{\partial}{\partial t}[\frac{1}
{\Gamma(1-\varepsilon)}f(t)]\pm\frac{\partial}{\partial t}[\alpha
\frac{\varepsilon(t)f(t)}{\Gamma(1-\varepsilon(t))}]
\end{equation}
were $\alpha$ defined by selection of regularization .The selection of a
sign in (\ref{eq6}) is defined  by a selection of the regularization. From
(\ref{eq10}) the opportunity follows (at the small fractional additives to
FD of time) to use for describing of changes of functions defined on
multifractal sets of time by means of using the renormalized ordinary
derivatives. At the same time, the dependence FD of the  time from
coordinates and time is concerned. Let's consider fractional dimension $d$
for case when $d$ smaller of unity ($d=1+\varepsilon$,$d<1$). For this
case fractional derivative (see (\ref{eq3}) for $n=1$) looks like
\begin{equation}\label{eq11}
D_{+,t}^{1-\varepsilon}f(t)=\frac{\partial}{\partial t}\int\limits_0^t
{\frac{f(\tau)d\tau}{\Gamma(\varepsilon(\tau)) (t-\tau\pm i\xi)^
{1-\varepsilon(\tau)}}}
\end{equation}
Taking into account, that for (\ref{eq11}) selection $\varepsilon$, by
virtue of definition (\ref{eq3})-(\ref{eq4})), is prohibited, for
including in (\ref{eq11}) value $D_{+,t}^{1-\varepsilon}f(t)$ at
$\varepsilon=0$ before a right member in (\ref{eq11}) (applicable only for
$\varepsilon>0$) it is necessary to take into account a addendum from
(\ref{eq9}) with $\varepsilon=0$, i.e. $\frac{\partial f(t)}{\partial t}$.
We receive (if use a rule of a regularization that was had used before for
$d>1$ i.e. $Reg(1/x)\to \alpha\delta(x)$ and relation
$\Gamma(1+\varepsilon)=\varepsilon\Gamma(\varepsilon))$
\begin{equation}\label{eq12}
\frac{\partial^{1-\varepsilon}}{\partial t^{1-\varepsilon}}f(t)\approx
\frac{\partial}{\partial t}f(t)\mp\frac{\partial}{\partial t}[\alpha
\frac{\varepsilon (r(t),t)f(t)}{\Gamma(1+\varepsilon (r(t),t))}f(t)]
\end{equation}
The approximate representation GFD by ordinary derivatives
(relations (\ref{eq8}),(\ref{eq10}),(\ref{eq12})) if use different
methods are very similar, so any of them may be used in follow
calculations. The above mentioned approximate connections of
generalized fractional derivatives (\ref{eq3})-(\ref{eq4}) defined
on the multifractal sets with fractional dimension $d({\bf
r}(t),t)$ (if $d({\bf r}(t),t)$ poorly distinguished from unity)
with ordinary derivatives may be ex-tend for the cases with
arbitrary $n$:$d({\bf r}(t),t)=n+\varepsilon({\bf r}(t),t)$,
$|\varepsilon|\ll 1$. Above mentioned reasoning make possible to
show, for a cases $\varepsilon\sim 1$ (but not close to integer
values), that the representations of generalized fractional
derivative by means of derivatives of integer order will contain
integer derivatives of arbitrary high orders. Let's consider a
symmetrical generalized fractional derivatives $D_{-,t}^{d} f(t)$
and $D_{+,t}^{d} f(t)$:
\begin{equation}\label{eq13}
D_{t}^{d}f(t)=0.5(D_{+,t}^{d}+D_{-,t}^{d})f(t)
\end{equation}
The symmetry of GFD allows to take into account the influence on
event that happens in the given instant  featured by function
$f(t)$ both past, and future (by fractional integration and
differentiation on time). For fractional integration and
differentiation at coordinates the symmetrical GFD takes into
account influence the event with given coordinate of all points of
space $$D_{{\bf r}}^d f(t)=0.5(D_{+,{\bf r}}^d+D_{-,{\bf
r}}^d)f(t)$$

At small difference of dimensions of time (or space) from unity
$D_{+,t}^{1+\varepsilon}\approx D_{-,t}^{1+\varepsilon}$ and so
on.

\section{Connection with covariant derivatives}

Let define (\ref{eq10}) as
\begin{equation}\label{eq14}
D_{+,t}^{1+\varepsilon}f(t)\approx A\frac{\partial}{\partial t}f-Bf
\end{equation}
where
\begin{equation}\label{eq15}
A({\bf r}(t),t)=\Gamma(1-\varepsilon)^{-1}+a\varepsilon
\end{equation}
\begin{equation}\label{eq16}
B({\bf r}(t),t)=
\Gamma(1-\varepsilon)^{-2}(1+a\varepsilon)\frac{\partial
\Gamma}{\partial t}-a\frac{\partial\varepsilon}{\partial
t},\quad(a=\pm 1)
\end{equation}
The relation (\ref{eq16}) reminds the covariant derivatives, frequently
meeting in physics. It is possible to show, that at a various selection of
a mathematical nature of function $f$ (vector, tensor etc.) and relevant
selection of function $\varepsilon$, GFD (\ref{eq10}) really coincides
with covariant derivatives (see \cite{15}, \cite{16}).

\section{Equations with generalized fractional derivatives}

The equations with GFD are possible to connect with natural
sciences (in particular, physics) when the fractal dimensions
$d_{t}$ and $d_{{\bf r}}$ are connected with describing
multifractal structure of a surfaces of solid bodies, structure of
chaos, structure of time and space (see, for example,
\cite{5},\cite{7},\cite{8},\cite{11},\cite{16}. In some cases GFD
are related to equations with FD that depends from functions (or
functionals) the same to which GFD was applied. It gives in the
interesting nonlinear fractional integral-differential equations
with GFD
\begin{equation}\label{eq17}
F(D_{+,t}^{d_t(f(t))})f(t)=0
\end{equation}
where $F$- function or functional from GFD.  A new class of the
equations in fractional integral-differential functionals
represent the equations such as (\ref{eq13}). Their examination,
apparently, is an interesting problem and represents a new
approach to describe problems of chaos.

\section{Conclusion}

The generalized Riemann-Liouville fractional  derivatives defined
in the paper allow to describe dynamics and changes of functions
defined on multifractal sets, in which every element of sets is
characterized by its own fractional dimension (depending on
coordinates and time). At small differences of fractional
dimensions from topological dimensions, generalized fractional
derivatives are represented through expressions similar to
covariant derivatives used in physics.


\begin{thebibliography}{99}
\bibitem{1} S.G.Samko, A.A.Kilbas , O.I.Marichev, \emph{Fractional
Integrals and Derivatives - Theory and Applications} (Gordon and Breach,
New York, 1993)
\bibitem{2} Schwarts L., \emph{Theorie des distribusions}, (Paris, Hermann, 1950,
vol.1,162p.; vol.2, 1951, 169p)
\bibitem{3} I.M.Gelfand, G.E.Shilov, \emph{Generalized Functions}
(Academic Press, New York, 1964)
\bibitem{4} Mandelbrot B. \emph{ The fractal geometry of nature} (W.H. Freeman,
New York, 1982)
\bibitem{5} L.Ya. Kobelev et al. \emph{Fractal Diffusion to a Rotating Disk}
Physics Doklady RAS, 1998, Vol.43, No.9, p.537
\bibitem{6}  Yu.L.Klimontovich  \emph{Statistical Theory of Open Systems}
(Kluwer, Dordrecht, 1995)
\bibitem{7} V.L.Kobelev et al. Physics  Doklady, 1998, Vol.43, No.8, p.487
\bibitem{8} L.Ya.Kobelev, et al, Physics Doklady, 1999, Vol.44, No.12
\bibitem{9} R.Metzler, E.Barkay, J.Klafter, Phys. Rev. Lett., v.82 (18),
p.3563
\bibitem{10} G.M.Zaslavsky, Chaos, 4,25, (1994)
\bibitem{11} Hausdorf F., Math. Ann. 79 (1919), P.157-179
\bibitem{12} Renyi A. \emph{Introduction to information theory, Appendix in:
Probability theory} (North Holland, Amsterdam, 1988)
\bibitem{13} Rudolph O., Fortshritte der Physik, v.43 (1995), No.5, P.349-450
\bibitem{14} Hadamard J., J.math.pures et appl.,Ser.4, 1892,T.8,P.101-186
\bibitem{15}Kobelev L.Ya.\emph{Fractal theory of time and space}(Ekaterinburg,
Konross, 1999)(in Rus.) KobelevL.Ya. \emph{The fractal theory of
time and space} Urals State University, Dep. v VINITI. 22.01.99,
No.189-B99
\bibitem{16}Kobelev L.Ya. \emph{Multifractality of time and special theory of
relativity} Urals State University, Dep. v VINITI 19.08.99,
No.2677-B99
\bibitem{17} L.Ya.Kobelev \emph{Multifractal time, covariant derivatives and
gauge invariance} Urals State University, Dep. v VINITI 07.09.99,
No.2791-B99
\end{thebibliography}
\end{document}